\documentclass[10pt]{article}
\pdfoutput=1
\usepackage[utf8]{inputenc}
\usepackage{amsfonts,amsmath,amsthm,amssymb}
\usepackage{nicefrac}
\usepackage{enumitem}
\usepackage[colorlinks,allcolors=blue]{hyperref}
\usepackage{url}
\usepackage[longnamesfirst, authoryear]{natbib}
\usepackage{breakcites}

\usepackage[ruled,vlined,noresetcount]{algorithm2e}
\SetKw{kwParams}{Parameters:}

\usepackage[title]{appendix}
\usepackage[margin=1in]{geometry}
\usepackage{tcolorbox}
\newtcolorbox{bluebox}{colback=blue!5!white,colframe=blue!75!black}
\usepackage{float}
\usepackage{tikz}
\usetikzlibrary{automata, positioning, arrows}
\usetikzlibrary{decorations.pathreplacing}

\newtheorem{theorem}{Theorem}

\newtheorem{lemma}{Lemma}

\theoremstyle{remark}

\newtheorem{example}{Example}

\let\Pr\relax
\DeclareMathOperator\Pr{\mathbb{P}}

\DeclareMathOperator\E{\mathbb{E}}

\DeclareMathOperator\Geom{Geometric}

\newcommand{\bx}{{\mathbf  x}}
\newcommand{\bX}{{\mathbf  X}}
\newcommand{\bd}{{\mathbf  d}}
\newcommand{\bc}{{\mathbf  c}}

\newcommand{\calA}{{\mathcal A}}

\newcommand{\calP}{{\mathcal P}}

\newcommand{\calS}{{\mathcal S}}

\newcommand{\Real}[1]{ { {\mathbb R}^{#1} } }

\title{\textbf{Bicausal Optimal Transport for Markov Chains via Dynamic Programming}}
\author{
Vrettos Moulos
\\
University of California Berkeley \\
\href{mailto:vrettos@berkeley.edu}{vrettos@berkeley.edu}
}
\date{\today}

\begin{document}

\maketitle

\begin{abstract}
In this paper we study the bicausal optimal transport problem for Markov chains, an optimal transport formulation suitable for stochastic processes which takes into consideration the accumulation of information as time evolves.
Our analysis is based on a relation between the transport problem and the theory of Markov decision processes.
This way we are able to derive necessary and sufficient conditions for optimality in the transport problem, as well as an iterative algorithm, namely the value iteration, for the calculation of the transportation cost.
Additionally, we draw the connection with the classic theory on couplings for Markov chains,
and in particular with the notion of faithful couplings.
Finally, we illustrate how the transportation cost appears naturally in the study of concentration of measure for Markov chains.
\end{abstract}

\section{Introduction}

Let $S$ be a finite set equipped with the discrete topology, and let $\calS$ be the corresponding Borel $\sigma$-field which in this case is the set of all subsets of $S$.
Let $S^\infty$ be the countable infinite product space, equipped with the product topology, and let $\calS^\infty$ be the corresponding Borel $\sigma$-field.
Let $\bX = (X_k)_{k=0}^\infty, ~\bX' = (X_k')_{k=0}^\infty$ be two discrete time stochastic processes on the measurable space $(S^\infty, \calS^\infty)$, governed by the probability laws $\mu, \mu'$ respectively.
By a coupling of $\bX, \bX'$ we mean a pair of stochastic processes $(\hat{\bX}, \hat{\bX}')$, on the measurable space $(S^\infty \times S^\infty, \calS^\infty \otimes \calS^\infty)$, governed by a probability law $\pi$ such that
$\pi (\cdot, S^\infty) = \mu$, and $\pi (S^\infty, \cdot ) = \mu'$.
Denote the set of all couplings of $\mu, \mu'$ by
\[
\Pi (\mu, \mu') = \left\{\pi \in \calP (S^\infty \times S^\infty) :
\pi (\cdot, S^\infty) = \mu,
~
\pi (S^\infty, \cdot) = \mu'\right\},
\]
where $\calP (S^\infty \times S^\infty)$ denotes the set of all probability laws on $(S^\infty \times S^\infty, \calS^\infty \otimes \calS^\infty)$.

Let $\bc : S^\infty \times S^\infty \to [0, \infty]$ be a $\calS^\infty \otimes \calS^\infty$-measurable cost function,
which has the following additive form
\begin{equation}
\bc (\bx, \bx') = \sum_{k=0}^\infty \beta^k c (x_k, y_k),
\end{equation}
for some $\beta \in (0, 1]$, and some $c : S \times S \to [0, \infty)$.
In particular, it will be of special interest the case that the cost function $\bc$ is the metric
\begin{equation}\label{eqn:metric}
    \bd_\beta (\bx, \bx') = \sum_{k=0}^\infty \beta^k I \{x_k \neq x_k'\},
\end{equation}
which for $\beta \in (0, 1)$ induces the product topology on $S^\infty$.
In a typical optimal transport problem, see for instance the book of~\cite{villani09},
we are interested in finding a coupling $\pi$
which minimizes the cost function $\bc$
according to the following cost criterion
\begin{equation}\label{eqn:classic-transport}
W (\mu, \mu') =
\inf_{\pi \in \Pi (\mu, \mu')}\:
\int \bc (\bx, \bx') \pi (d \bx, d \bx').
\end{equation}

Such a formulation might be inadequate in the context of stochastic processes, as the evolution over time matters, and has to be accounted.
In the context of finite horizon processes
the works of~\cite{pflug09,pflug12} recognize this and introduce the nested distance which takes into consideration the filtrations. 
They are motivated by applications of the nested distance for scenario reduction in the context of multistage stochastic optimization.
See also the work of~\cite{backhoff20a} for an application of the nested distance to stability in mathematical finance.
The metric and topological properties of the nested distance has been recently studied in~\cite{backhoff20b}.
A generalization of the nested distance is the causal optimal transport problem introduced by~\cite{lassalle15} and further developed by~\cite{backhoff17} where a dynamic programming principle is developed as well.
All those works deal in great generality with processes of finite horizon.
In this paper we study the bicausal optimal transport problem for Markov chains over an infinite horizon, drawing motivation by the classic theory of couplings for Markov chains where one might naturally seek for an adapted coupling, i.e. one that cannot see into the future, that couples
two chains as fast as possible.

To introduce the bicausal optimal transport problem we first note that any probability law $\pi$ can be factorized, for every $n$, as
\begin{equation}\label{eqn:factorization}
\begin{split}
&\pi ((d x_0, \ldots, d x_n), (d x_0', \ldots, d x_n')) \\
&\qquad = 
\pi (d x_0, d x_0')
\pi (d x_1, d x_1' \mid x_0, x_0')
\cdots
\pi (d x_n, d x_n' \mid x_0, \ldots, x_{n-1}, x_0', \ldots, x_{n-1}'),
\end{split}
\end{equation}
where $\pi (\cdot \mid x_0, \ldots, x_{k-1}, x_0', \ldots, x_{k-1}')$ denotes the conditional probability law of $((\hat{X})_{i=k}^\infty, (\hat{X}')_{i=k}^\infty)$ given that
$\hat{X}_0 = x_0, \ldots, \hat{X}_{k-1} = x_{k-1}, \hat{X}_0' = x_0', \ldots, \hat{X}_{k-1}' = x_{k-1}'$. In this work we are interested in couplings $\pi$ which are bicausal, in the sense that for every $n, ~\pi (d x_n, S \mid x_0, \ldots, x_{n-1}, x_0', \ldots, x_{n-1}') = \mu (d x_n \mid x_0, \ldots, x_{n-1})$, and
$\pi (S, d x_n' \mid x_0, \ldots, x_{n-1}, x_0', \ldots, x_{n-1}') = \mu' (d x_n' \mid x_0', \ldots, x_{n-1}')$.
We denote the set of all bicausal couplings of $\mu, \mu'$ by
\[
\Pi_{bc} (\mu, \mu') = \left\{\pi \in \calP (S^\infty \times S^\infty) :
\begin{array}{l}
\pi (d x_n, S \mid x_0, \ldots, x_{n-1}, x_0', \ldots, x_{n-1}') = \mu (d x_n \mid x_0, \ldots, x_{n-1}) \\
\pi (S, d x_n' \mid x_0, \ldots, x_{n-1}, x_0', \ldots, x_{n-1}') = \mu' (d x_n' \mid x_0', \ldots, x_{n-1}')
\end{array}, ~\text{for every}~n
\right\}.
\]
Due to the factorization~\eqref{eqn:factorization} it is clear that $\Pi_{bc} (\mu, \mu') \subseteq \Pi (\mu, \mu')$.
Additionally, the product measure $\mu \otimes \mu' \in \Pi_{bc} (\mu, \mu')$, hence none of those sets is empty.
The corresponding bicausal optimal transport problem can be written as
\begin{equation}\label{eqn:causal-transport}
W_{bc} (\mu, \mu') =
\inf_{\pi \in \Pi_{bc} (\mu, \mu')}\:
\int \bc (\bx, \bx') \pi (d \bx, d \bx').
\end{equation}

The bicausal optimal transport problem~\eqref{eqn:causal-transport} is particular interesting in the case that $\mu, \mu'$ are Markovian laws, i.e. $(X_k)_{k=0}^\infty, (X_k')_{k=0}^\infty$ are Markov chains. For the rest of this paper we assume that there are initial states $x_0, x_0' \in S$, and transition probability kernels $P, P' : S \times \calS \to [0, 1]$ such that for every $n$
\begin{align*}
    \mu (\{x_0\}, d x_1, \ldots, d x_n) &=
    P (x_0, dx_1) P (x_1, dx_2) \cdots P (x_{n-1}, dx_n),
    \\
    \mu' (\{x_0'\}, d x_1', \ldots, d x_n') &=
    P' (x_0', dx_1') P' (x_1', dx_2') \cdots P' (x_{n-1}', dx_n'),
\end{align*}
and we write
\[
\mu = \mathrm{Markov} (x_0, P), \quad \mu' = \mathrm{Markov} (x_0', P').
\]
We note that using two fixed initial states $x_0, x_0'$ is as general as considering arbitrary initial distributions, since $x_0, x_0'$ can be thought of as auxiliary states inducing arbitrary initial distributions, $P (x_0, \cdot), P (x_0', \cdot)$ to $X_1, X_1'$ respectively.
For extra clarity we rewrite~\eqref{eqn:causal-transport} as
\begin{equation}\label{eqn:causal-transport-Markov}
W_{bc} \left(\mathrm{Markov} (x_0, P), \mathrm{Markov} (x_0', P')\right) =
\inf_{\pi \in \Pi_{bc} \left(\mathrm{Markov} (x_0, P), \mathrm{Markov} (x_0', P')\right)}\:
\int \bc (\bx, \bx') \pi (d \bx, d \bx').
\end{equation}
We study the transportation problem~\eqref{eqn:causal-transport-Markov} in~\autoref{sec:DP} under the lens of dynamic programming, where we develop optimality conditions, as well as an iterative procedure, namely the value iteration, that solves the transportation problem.

\paragraph{Motivation}

In the special case that $P = P'$, and the cost function $\bc$ is the metric
$\bd_\beta$ given in~\eqref{eqn:metric},
the transportation problem~\eqref{eqn:classic-transport} can be solved explicitly
by works on maximal coupling from~\cite{griffeath1975,pitman1976,goldstein1979}.
In particular
\[
W \left(\mathrm{Markov} (x_0, P), \mathrm{Markov} (x_0', P)\right) = 
\sum_{k=0}^\infty \beta^k \left\|
P^k (x_0, \cdot) - P^k (x_0', \cdot)
\right\|_{TV},
\]
where $\|\cdot\|_{TV}$ stands for half the total variation norm of a signed measure.

When $\beta = 1$ the transportation problem~\eqref{eqn:classic-transport} reduces to finding a coupling of two different initializations of the same Markov chain that couples them in the smallest expected time. 
By introducing the coupling time
\[
T = \inf\:\left\{n \ge 0 : \hat{X}_n = \hat{X}_n' \right\},
\]
the transportation problem \eqref{eqn:classic-transport} can be written as
\begin{equation}\label{eqn:transport-concentration}
W \left(\mathrm{Markov} (x_0, P), \mathrm{Markov} (x_0', P)\right)  = 
\inf_{\pi \in \Pi \left(\mathrm{Markov} (x_0, P), \mathrm{Markov} (x_0', P)\right)}\:
\E^\pi \left[
T
\right].
\end{equation}
This transportation problem is particularly important because it directly
leads to a bounded differences concentration inequality for Markov chains
as we discuss in~\autoref{sec:concentration}.

The maximal coupling of~\cite{pitman1976}, works on the space-time plane by first simulating the meeting point, and then constructing the forward and backward chains. As such the coupling `cheats' by looking into the future. \cite{rosenthal1997} initiates the discussion of faithful couplings that do not look into the future,
motivated by the fact that such couplings automatically posses the `now equals forever' property
which roughly speaking says that the two chains becoming equal at a single time is equivalent to having them remain equal for all future times.
It is the bicausal version of the transportation problem
\begin{equation}\label{eqn:casual-transport-Markov-coupling}
W_{bc} (\left(\mathrm{Markov} (x_0, P), \mathrm{Markov} (x_0', P')\right)) = 
\inf_{\pi \in \Pi_{bc} (\left(\mathrm{Markov} (x_0, P), \mathrm{Markov} (x_0', P')\right))}\:
\E^\pi \left[
T
\right],
\end{equation}
that seeks for faithful couplings, that do not look into the future, and minimize the expected coupling time. In~\autoref{sec:coupling-MC} we provide necessary and sufficient conditions for optimality at the transportation problem~\eqref{eqn:casual-transport-Markov-coupling}, as well as a discussion about properties of optimal couplings.

\section{Markovian Couplings}

Among the set of all bicausal couplings, $\Pi_{bc} \left(\mathrm{Markov} (x_0, P), ~\mathrm{Markov} (x_0', P')\right)$, it suffices to turn our attention to Markovian couplings
when considering the transportation problem~\eqref{eqn:causal-transport}. We will establish this in~\autoref{sec:DP} as a consequence of the dynamic programming theory. A Markovian coupling of $\mathrm{Markov} (x_0, P), ~\mathrm{Markov} (x_0', P')$ is specified by a transition kernel $Q : (S \times S) \times (\calS \otimes \calS) \to [0, 1]$ such that
$Q ((x, x'), (\cdot, S)) = P (x, \cdot)$, and
$Q ((x, x'), (S, \cdot)) = P' (x', \cdot)$. We denote the set of all such transition kernels by
\[
\Pi_M (P, P') = \left\{
Q : 
Q ((x, x'), (\cdot, S)) = P (x, \cdot),~
Q ((x, x'), (S, \cdot)) = P' (x', \cdot)
\right\}
\]
The corresponding Markovian coupling is given, for every $n$, by
\begin{equation}\label{eqn:markovian-coupling}
\pi ((\{ x_0\}, \ldots, d x_n), (\{x_0'\}, \ldots, d x_n')) =
Q ((x_0, x_0'), (d x_1, d x_1')) \cdots Q ((x_{n-1}, x_{n-1}'), (d x_n, d x_n')).
\end{equation}
We note that any Markovian coupling is bicausal.

We now present some basic examples of Markovian couplings,
for the case of a single Markov chain $P = P'$,
which have been used extensively in the coupling literature, see for instance the book~\cite{thorisson00}.
\begin{example}
The classic coupling, initially introduced by Doeblin in order to establish the convergence theorem for Markov chains, asserts that $\hat{X}_n$ and $\hat{X}_n'$ evolve independently until they reach a common state for the first time, and afterwards they move identically.
\begin{equation}\label{eqn:sticky-coupling}
Q_{\mathrm{classic}} ((x, x'), (y, y')) =
\begin{cases}
P (x, y) P (x', y'), & \text{if}~ x \neq x', \\
P (x, y), & \text{if}~ x = x', ~\text{and}~ y = y', \\
0, & \text{otherwise},
\end{cases}
\end{equation}
\end{example}
\begin{example}
A variant of the classic coupling asserts that $\hat{X}_n$ and $\hat{X}_n'$ evolve independently at all times,
and this independent coupling can be used as well to establish the convergence theorem for Markov chains. 
\begin{equation}\label{eqn:ind-coupling}
Q_{\mathrm{ind}} ((x, x'), (y, y')) =
P (x, y) P (x', y'),
\end{equation}
\end{example}
In both the classic and the independent coupling it is apparent that if we seek for a Markovian coupling that minimizes some cost criterion, e.g. attaining coupling at the smallest expected time, then the independent movement can be wasteful. Instead, one should coordinate the movement of the two copies in a way that optimizes the objective under consideration. A first such attempt is the following coupling attributed to Wasserstein. 
\begin{example}
Given that $\hat{X}_{n-1} = x$ and $\hat{X}_{n-1}' = x'$, the Wassertstein coupling makes $\hat{X}_n$
and $\hat{X}_n'$ agree with as great probability as possible (which is $1 - \|P (x, \cdot) - P (x', \cdot)\|_{TV}$), and then given that they differ they are made conditionally independent.
\begin{equation}\label{eqn:w-coupling}
    Q_{\mathrm{W}} ((x, x'), (y, y')) =
    \begin{cases}
    0, & \text{if} ~ x = x', ~\text{and}~ y \neq y', \\
    P (x, y), & \text{if} ~ x = x', ~\text{and}~ y = y', \\
    P (x, y) \land P (x', y), & \text{if} ~ x \neq x', ~\text{and}~ y = y', \\
    0, & \text{if} ~ x \neq x', ~ y \neq y', ~\text{and}~ \|P (x, \cdot) - P (x', \cdot)\|_{TV} = 0,
    \\
    \frac{\left(P (x, y) - P (x', y)\right)^+
    \left(P (x', y') - P (x, y')\right)^+}{\|P (x, \cdot) - P (x', \cdot)\|_{TV}}, & \text{if} ~ x \neq x', ~ y \neq y', ~\text{and}~ \|P (x, \cdot) - P (x', \cdot)\|_{TV} \neq 0,
    \end{cases}
\end{equation}
where $a \land b = \min (a, b), ~ a^+ = - ((-a) \land 0),$ and $a^- = a \land 0$.
\end{example}
Making $\hat{X}_n$ and $\hat{X}_n'$ agree with as great probability as possible is a good first step towards a Markovian coupling with the smallest expected coupling time, although it might be too greedy of a choice and the conditional independence assertion surely leaves more room for improvements. It is the objective of this paper to provide a characterization of optimal Markovian couplings using the theory of dynamic programming.
\section{Bicausal Optimal Transport for Markov Chains via Dynamic Programming}\label{sec:DP}

We start by illustrating that the bicasual transport problem for Markov chains~\eqref{eqn:causal-transport-Markov} can be viewed as an instance of infinite horizon dynamic programming.
When $\beta \in (0, 1)$ we have an instance of discounted dynamic programming, initially studied by~\cite{blackwell1965}, while when $\beta = 1$ we have an instance of negative dynamic programming, initially studied by~\cite{strauch1966}.
For two Markov chains $\mathrm{Markov} (x_0, P), ~ \mathrm{Markov} (x_0', P')$ on the same state space $(S, \calS)$, and for the bicausal optimal transport problem~\eqref{eqn:causal-transport-Markov} between them, the associated underlying Markov decision process is given by the tuple 
$((S \times S, \calS \otimes \calS), (A, \calA), U, q, \beta, c)$ where:
\begin{itemize}
    \item $(S \times S, \calS \otimes \calS)$ stands for the state space of the Markov decision process.
    
    \item $(A, \calA)$ is the action space, where $A = \calP (S \times S)$ is the set of all probability distributions on $S \times S$ equipped with the subspace topology induced from the standard topology on $\Real{|S \times S|}$, and $\calA$ stands for the corresponding Borel $\sigma$-field.
    
    \item $U (x,x') = \left\{a \in A : a (\cdot, S) = P (x, \cdot), ~ a (S, \cdot) = P' (x', \cdot)\right\}$ is the set of all allowable actions at state $(x, x')$,
    i.e. all the probability distributions on $S \times S$ which respect the coupling constraints.
    
    \item $q (\cdot \mid (x, x'), a) = a$ is the distribution of the state next visited by the Markov decision process if the system is currently in state $(x, x')$ and action $a \in U (x, x')$ is taken.
    
    \item $\beta \in (0, 1]$ is the discount factor.
    
    \item $c : S \times S \to [0, \infty)$ is the cost function.
    
\end{itemize}
A policy is a bicausal coupling $\mu \in \Pi_{bc} (\mathrm{Markov} (x_0, P), \mathrm{Markov} (x_0', P'))$,
and it can be decomposed as a sequence of conditional distributions as in~\eqref{eqn:factorization} so that if the coupling $\mu$ is used,
and up to time $n$ we observe the trajectory
$x_0, \ldots, x_n, x_0', \ldots, x_n'$ then the action taken at time $n$ is $\mu (\cdot \mid x_0, \ldots, x_n, x_0', \ldots, x_n')$ which is also the distribution of the state visited by the Markov decision process at time $n+1$. As it turns out there exists an optimal coupling for which the conditional distributions do not depend on the whole trajectory but just on the current state,
i.e. there exists an optimal coupling which is Markovian in the sense of~\eqref{eqn:markovian-coupling}.

We proceed with the definition of some typical operators from the dynamic programming literature.
Let $F$ be the set of all extended real valued functions $V : S \times S \to [0, \infty]$.
For $Q \in \Pi_M (P, P')$ define the operator $T_Q : F \to F$ by
\[
T_Q (V) (x, x') = c (x, x') + 
\beta \int Q \left((x, x'), (dy, dy')\right) V (y, y'),
\]
which we may also write in functional notation as
\[
T_Q (V) = c + \beta Q V.
\]
Additionally, define the Bellman operator $T : F \to F$ by
\[
T (V) (x, x') = c (x, x') +
\beta \inf_{a \in U (x, x')}\: \int a (dy, dy') V (y, y'),
\]
or in functional notation as
\[
T (V) = c + \beta \inf_{Q \in \Pi_M (P, P')}\: Q V.
\]
We note that when $\beta \in (0, 1)$ the Bellman operator $T$ is a $\beta$-contraction with respect to the $\sup$-norm on $F$,
and when $\beta = 1$ the Bellman operator $T$ is an increasing mapping. In the following theorem we summarize the main results from this dynamic programming interpretation of the bicausal optimal transport problem between two Markov chains~\eqref{eqn:causal-transport-Markov}. As it is typical in dynamic programming we consider the transportation cost $W_{bc} \left(\mathrm{Markov} (x_0, P), \mathrm{Markov} (x_0', P')\right)$ as a function
$W_{bc} : S \times S \to [0, \infty]$ of the initializations of the two Markov chains,
and we write $W_{bc} (x_0, x_0')$ for the optimal cost.

\begin{theorem}\label{thm:DP}\hfill\break
\begin{enumerate}
    \item The transportation cost $W_{bc}$ is a solution to the fixed point equation $V = T (V)$. 
    When $\beta \in (0, 1)$ it is the unique solution,
    while when $\beta = 1$ if $V \ge 0$ and $V = T (V)$ then $V \ge W_{bc}$.
    
    \item The transportation cost $W_{bc}$ can be calculated via the fixed point iteration
    \begin{equation}\label{eqn:value-iteration}
    \begin{cases}
    V_0 &= 0, \\
    V_k &= T (V_{k-1}), \quad k = 1, 2, \ldots,
    \end{cases}
    \end{equation}
    If $\beta \in (0, 1)$, then $\|V_k - W_{bc}\|_\infty \le \beta^k \|W_{bc}\|_\infty$, and thus $V_k \to W_{bc}$ as $k \to \infty$ at a linear rate.
    If $\beta = 1$, then the convergence is monotonic, $V_k \uparrow W_{bc}$ as $k \to \infty$.
    
    \item There exists an optimal Markovian coupling.
    
    \item $Q$ is an optimal Markovian coupling if and only if $T_Q (W_{bc}) = W_{bc}$.
\end{enumerate}
\begin{proof}
When $\beta \in (0, 1)$ parts 1, 2 and 4 follow from Proposition 1 in~\cite{bertsekas1977}.
When $\beta = 1$ parts 1 and 3 follows from Propositions 5 and 7 in~\cite{bertsekas1977}.

For the rest we need to note that for every $x, x' \in S, ~ \lambda \in [0, \infty)$ and $k$, the set
\[
U_k ((x, x'), \lambda) = 
\left\{
a \in U (x, x') : 
c (x, x') + \beta \int a (dy, dy') V_k (y, y') \le \lambda
\right\},
\]
is compact as the intersection of the compact set $U (x, x')$, with a closed half-space.

Then for $\beta \in (0, 1)$ part 3 follows from Proposition 14 in~\cite{bertsekas1977},
and for $\beta = 1$ parts 2 and 3 follow from Proposition 12 in~\cite{bertsekas1977}.
\end{proof}
\end{theorem}

We note that due to the special structure of $U (x, x')$, 
it is a convex polytope arising from the intersection of a probability simplex with an affine space,
the value iteration~\eqref{eqn:value-iteration} proceeds by solving at each iteration a linear program.
Thus the value iteration~\eqref{eqn:value-iteration} in this case can be thought as sequence of finite dimensional linear programs approximating the bicausal optimal transport cost $W_{bc}$ which in~\eqref{eqn:causal-transport-Markov} is formulated as an infinite dimensional linear program.
\section{Coupling of Markov Chains in Minimum Expected Time}\label{sec:coupling-MC}

In this section we specialize the bicausal optimal transport for Markov chains to the case that we have a single irreducible and aperiodic transition kernel $P$, and the cost function $\bc$ is the metric $\bd_1$. So essentially we are solving for the bicausal coupling that couples two Markov chains with different initializations and the same transition kernel in the smallest expected time
\begin{align*}
W_{bc} (\mathrm{Markov} (x_0, P), \mathrm{Markov} (x_0', P)) 
&=
\inf_{\pi \in \Pi_{bc} \left(\mathrm{Markov} (x_0, P), \mathrm{Markov} (x_0', P')\right)}\:
\int \bd_1 (\bx, \bx') \pi (d \bx, d \bx') \\
&=
\inf_{\pi \in \Pi_{bc} (\left(\mathrm{Markov} (x_0, P), \mathrm{Markov} (x_0', P')\right))}\:
\E^\pi \left[
T
\right].
\end{align*}
Although in the general framework of negative dynamic programming the fixed point equation
$V = T (V)$ is only a necessary condition for optimality,
in our specialized setting we can establish that it is also sufficient,
giving thus a complete set of necessary and sufficient conditions for both the optimal cost $W_{bc}$, and the optimal Markovian coupling $Q$.
\begin{theorem}
$W_{bc}$ is the unique solution of the equations
\begin{equation}\label{eqn:optimality}
    0 \le V < \infty, \quad
    V = T (V), \quad
    \text{and} \quad
    V (x, x) = 0 \quad
    \text{for $x \in S$}.
\end{equation}
\begin{proof}
We already know from~\autoref{thm:DP} that $W_{bc} = T (W_{bc})$.
Using the classic Markovian coupling,~\eqref{eqn:sticky-coupling}
we see that $W_{bc} (x, x) = 0$ for all $x \in S$.
Using the independent Markovian coupling,~\eqref{eqn:ind-coupling},
which induces an irreducible Markov chain on $S \times S$ we see that
\[
W_{bc} (x, x') \le
\min_{y \in S}\: \E_{(x, x')}^{Q_{\mathrm{ind}}} [T_{(y,y)}]
< \infty,
\quad \text{for all}~ x,x' \in S.
\]

Let $V : S \to [0, \infty)$ be any function satisfying equations~\eqref{eqn:optimality}.
Let $Q$ be a Markovian coupling such that
$W_{bc} = c + Q W_{bc}$. Then
\begin{equation}\label{eqn:basic-ineq}
V = T (V) \le
c + Q V =
c + Q T (V) \le
c + Q c + Q^2 V \le 
\cdots \le
\sum_{k=0}^{n-1} Q^k c + Q^n V.
\end{equation}
By definition of $Q$ we have that
\begin{equation}\label{eqn:opt-cost}
W_{bc} = \sum_{k=0}^\infty Q^k c,
\end{equation}
and since $W_{bc} < \infty$ we see that $\lim_{n \to \infty} Q^n c = 0$.
We clearly have that $c \le V$, and because $c (x, x') = 0 \;\Rightarrow\; V (x, x') = 0$, we obtain that $c \le V \le \|V\|_\infty c$. 
Since $V < \infty$ we deduce that
\begin{equation}\label{eqn:final-cost}
\lim_{n \to \infty} Q^n V = 0.
\end{equation}
Combining~\eqref{eqn:basic-ineq},~\eqref{eqn:opt-cost}, and~\eqref{eqn:final-cost} we obtain that $V \le W_{bc}$, and because from~\autoref{thm:DP} $W_{bc}$ is the minimal fixed point we deduce that $V = W_{bc}$.
\end{proof}
\end{theorem}
Next we dig in some properties of an optimal Markovian coupling.
In particular, we show that an optimal Markovian coupling enjoys the `sticky' property of the classic coupling~\eqref{eqn:sticky-coupling}, i.e. under an optimal Markovian coupling the two chains evolve in the same way as soon as they meet.
\begin{lemma}\label{lem:sticky}
Any optimal Markovian coupling $Q$
sticks to the diagonal as soon as it hits it,
i.e.
\[
Q ((x, x), (y, y')) = I \{y = y'\} P (x, y).
\]
\begin{proof}
Fix an optimal Markovian coupling $Q$. From~~\autoref{thm:DP} it satisfies the equation
\[
c + Q W_{bc} = W_{bc},
\]
and so in particular
\[
\int Q ((x, x), (d y, d y')) W_{bc} (y, y') = 0.
\]
Since for $y \neq y', ~ W_{bc} (y, y') \ge 1$ we have that $Q ((x, x), (y, y')) = 0$.
Then it follows from the coupling constraint that $Q ((x, x), (y, y)) = P (x, y)$.
\end{proof}
\end{lemma}
Additionally, we show that for two state chains the Wasserstein coupling~\eqref{eqn:w-coupling} is the only optimal Markovian coupling.
\begin{lemma}
When $|S| = 2$ there is a unique optimal Markovian coupling which is precisely the Wasserstein coupling~\eqref{eqn:w-coupling}.
\begin{proof}
Let $S = \{x, x'\}$.
Due to symmetry we have that that $W_{bc} (x, x') = W_{bc} (x', x)$,
and in addition $W_{bc} (x, x) = W_{bc} (x', x') = 0$.
So from~\autoref{thm:DP} we get that
\begin{equation}\label{eqn:two-state}
W_{bc} (x, x') = 
1 +
\min_{a \in U (x, x')}\:
(1 - a (x, x) - a (x', x')) W_{bc} (x, x').
\end{equation}
It is clear that in the minimization in~\eqref{eqn:two-state} we need to make 
$a (x, x)$, and $a (x', x')$ as large as possible.
Due to the coupling constraint, $a \in U (x, x')$, those largest values are
\[
a (x, x) = P (x, x) \land P (x', x),
\quad\text{and}\quad
a (x', x') = P (x, x') \land P (x', x').
\]
Then from the coupling constraints there are unique corresponding values for $a (x, x')$, and $a (x', x)$. In particular
\begin{align*}
a (x, x') &= 
(P (x,x) - P (x', x))^+ = (P (x', x') - P (x, x'))^+, \\
a (x', x) &=
(P (x', x) - P (x, x))^+ = (P (x, x') - P (x', x'))^+.
\end{align*}
We further note that
\[
\|P (x, \cdot) - P (x', \cdot)\|_{TV} = |P (x, x) - P (x', x)| = |P (x', x') - P (x, x')|.
\]
Hence in conjunction with~\autoref{lem:sticky} we conclude that there exists a unique optimal Markovian coupling and this is the Wassertstein coupling~\eqref{eqn:w-coupling}.

Moreover, we have closed form expressions for both the non-causal and the bicausal optimal transport costs
\begin{align*}
& W (x, x') =
\frac{|P (x, x') - P (x', x)|}{P (x, x') + P (x', x)} \cdot 
\frac{1}{1 - \|P (x, \cdot) - P (x', \cdot)\|_{TV}} \\
& \quad <
W_{bc} (x, x') = 
\frac{1}{1 - \|P (x, \cdot) - P (x', \cdot)\|_{TV}}.
\end{align*}
\end{proof}
\end{lemma}
Finally, we give an easy upper bound on the bicausal optimal transport cost for contractive Markov chains.
\begin{lemma}\label{lem:geometric}
Let $\delta (P) = \max_{x, x' \in S}\: \|P (x, \cdot) - P (x', \cdot)\|_{TV}$ be the Doeblin-Dobrushin contraction coefficient, and assume that $\delta (P) < 1$. Then
\[
\|W_{bc}\|_\infty \le \frac{1}{1-\delta (P)}.
\]
\begin{proof}
By definition $W_{bc} (x_0, x_0')$ is upper bounded by the cost incurred when the Wasserstein coupling~\eqref{eqn:w-coupling} is used.
Under the Wasserstein coupling, $Q_W$, the probability that we hit the diagonal in one step from state
$(x, x')$ is $1-\|P  (x, \cdot) - P (x', \cdot)\|_{TV}$.
Thus $\bd_1 (\hat{\bX}, \hat{\bX}')$ under the Markovian coupling induced by $Q_W$ is stochastically dominated by $\Geom (1-\delta (P))$, and thus
\[
W_{bc} (x_0, x_0') \le 
\E_{(x_0, x_0')}^{Q_W} [\bd_1 (\hat{\bX}, \hat{\bX}')] \le 
\frac{1}{1-\delta (P)},
\quad \text{for any}~ x_0, x_0' \in S.
\]
\end{proof}
\end{lemma}
\section{Concentration of Measure for Markov Chains}\label{sec:concentration}

In this section we demonstrate how the transportation cost~\eqref{eqn:transport-concentration}
shows up naturally when one studies concentration of measure for Markov chains.
This was first observed in the works of~\cite{Marton-96-a,Marton-96-b} about contracting Markov chains,
and thereafter greatly generalized for classes of dependent random processes in terms of various mixing coefficients using the transportation-information method~\cite{Marton-98, Marton-03, Samson-00}.
\cite{chazottes07} uses the martingale method combined with maximal coupling to derive concentration for dependent processes, while~\cite{kontorovich08} uses the martingale method and a linear programming inequality for the same task.
For the Markovian case~\cite{Paulin15} using the martingale method establishes a concentration inequality that involves the mixing time of the chain.

Let $f : S^n \to \Real{}$ be a functions which is $1$-Lipschitz with respect to the Hamming distance
\[
f (x_1, \ldots, x_n) - f (x_1', \ldots, x_n') \le \sum_{k=1}^n I \{x_k \neq x_k'\},
\quad \text{for}~ x_1, \ldots, x_n, x_1', \ldots, x_n' \in S.
\]
Let $\bX \sim \mathrm{Markov} (x_0, P)$, where $P$ is an irreducible and aperiodic transition kernel. We would like to study the deviation of the random variable
$f (X_1, \ldots, X_n)$ from its mean $\E_{x_0}^P [f (X_1, \ldots, X_n)]$.
The typical approach to do so is the martingale method.
For $k = 0, \ldots, n$ we define the martingale
\[
Z_k = \E_{x_0}^P \left( f (X_1, \ldots, X_n) \middle| X_1, \ldots, X_k \right),
\]
and for $k=1, \ldots, n$ we define the martingale differences
\[
\Delta_k = Z_k - Z_{k-1}.
\]
Then the quantity of which we want to control the deviations can be written as a telescoping sum of the martingale differences
\[
f (X_1, \ldots, X_n) - \E_{x_0}^P\left[f (X_1, \ldots, X_n)\right] = \sum_{k=1}^n \Delta_k,
\]
and it suffices to control the range of the martingale differences.
For this we note that
\[
\min_{x' \in S}\: \left\{
\E_{x_0}^P \left( f (X_1, \ldots, X_n) \middle| X_1, \ldots, X_k = x' \right)
\right\}
- \E_{x_0}^P \left( f (X_1, \ldots, X_n) \middle| X_1, \ldots, X_{k-1} \right)
\le \Delta_k,
\]
and that
\[
 \Delta_k \le
\max_{x \in S}\: \left\{
\E_{x_0}^P \left( f (X_1, \ldots, X_n) \middle| X_1, \ldots, X_k = x \right)
\right\}
- \E_{x_0}^P \left( f (X_1, \ldots, X_n) \middle| X_1, \ldots, X_{k-1} \right).
\]
Thus in order to bound the length of the range of the martingale difference, we just need to bound
\begin{equation}\label{eqn:variation}
\max_{x,x' \in S}\: \left\{
\E_{x_0}^P \left( f (X_1, \ldots, X_n) \middle| X_1, \ldots, X_k = x \right) - 
\E_{x_0}^P \left( f (X_1, \ldots, X_n) \middle| X_1, \ldots, X_k = x' \right)
\right\}.
\end{equation}
Fix $x_1, \ldots, x_{k-1}, x, x' \in S$, and a coupling
$\pi \in \Pi (\mu (\cdot \mid x_1, \ldots, x_{k-1}, x), \mu (\cdot \mid x_1, \ldots, x_{k-1}, x'))$,
where $\mu = \mathrm{Markov} (x_0, P)$. Then
\begin{align*}
&\E_{x_0}^P \left( f (X_1, \ldots, X_n) \middle| X_1 = x_1, \ldots, X_{k-1} = x_{k-1}, X_k = x \right) - 
\E_{x_0}^P \left( f (X_1, \ldots, X_n) \middle| X_1 = x_1, \ldots, X_{k-1} = x_{k-1}, X_k = x' \right) \\
&\quad =
\E^\pi \left[
f (x_1, \ldots, x_{k-1}, \hat{X}_k, \ldots, \hat{X}_n) -
f (x_1, \ldots, x_{k-1}, \hat{X}_k', \ldots, \hat{X}_n')
\right] \\
&\quad\le 
\E^\pi \left[
\sum_{i=k}^n I \{\hat{X}_i \neq \hat{X}_i'\}
\right] \le 
\E^\pi \left[
\sum_{i=k}^\infty I \{\hat{X}_i \neq \hat{X}_i'\}
\right],
\end{align*}
where we used the Lipschitz condition for $f$.
Thus minimizing over the coupling $\pi$ we obtain that
\[
\E_{x_0}^P \left( f (X_1, \ldots, X_n) \middle| X_1, \ldots, X_k = x \right) - 
\E_{x_0}^P \left( f (X_1, \ldots, X_n) \middle| X_1, \ldots, X_k = x' \right) 
\le W (x, x').
\]
All in all, we can bound the length of the range of the martingale difference $\Delta_k$ by
$\|W\|_\infty$. Then using the standard martingale method one can obtain the concentration inequality
\begin{equation}\label{eqn:concentration-ineq}
\Pr_{x_0}^P \left(\left|
f (X_1, \ldots, X_n) - \E_{x_0}^P [f (X_1, \ldots, X_n)]
\right|\ge t\right) \le 2 \exp\left\{-\frac{2 t^2}{n \|W\|_\infty^2}\right\}.
\end{equation}
For a full derivation of a concentration inequality which works in more general dependent settings, than just Markovian dependence, the interested reader is refereed to Theorem 1 in~\cite{chazottes07} which uses the martingale method together with maximal coupling, and to
Theorem 1.1 of~\cite{kontorovich08} which uses a linear programming inequality instead of a coupling argument. Additionally, we note that if the Markov chain is periodic and thus coupling techniques are not any more applicable one can still bound~\eqref{eqn:variation}, in the special case that the function $f$ is additive, by using hitting time arguments as it is done in~\cite{moulos20}.

Clearly $\|W\|_\infty \le \|W_{bc}\|_\infty$, and thus replacing $\|W\|_\infty$ with $\|W_{bc}\|_\infty$ in~\eqref{eqn:concentration-ineq} results in a weaker inequality,
although in this way the variance proxy $\|W_{bc}\|_\infty^2$ has the following interpretation: let $Q$ be an optimal Markovian coupling, then $W_{bc} (x, x')$ corresponds to the expected time to hit diagonal when we start from $(x, x')$ and we transition according $Q$, thus the variance proxy is the squared expected time required to hit the diagonal when the least favorite initialization is used. Additionally, when the transition kernel is contracting, $\delta (P) < 1$, we can apply~\autoref{lem:geometric} and further replace $\|W_{bc}\|_\infty$ by $1/(1-\delta (P))$, which results in a specialized version of Theorem 1.2 in~\cite{kontorovich08}.

\section*{Acknowledgements}
We would like to thank Venkat Anantharam, Sinho Chewi, and Satish Rao for helpful discussions.
This research was supported in part by the NSF grant CCF-1816861

\bibliographystyle{apalike}
\bibliography{references}

\end{document}